\documentclass[11pt, twoside]{article}
\usepackage{amsmath,amsthm}
\usepackage{amsfonts}
\usepackage{amssymb,latexsym}
\usepackage{enumerate}
\newtheorem{theorem}{Theorem}[section]
\newtheorem{remark}[theorem]{Remark}
\newtheorem{proposition}[theorem]{Proposition}
\newtheorem{lemma}[theorem]{Lemma}

\newtheorem{corollary}[theorem]{Corollary}

\def\ind{1{\hskip -3 pt}\hbox{\textsc{I}}}
\def\n{\noindent}
\def\n{\noindent}

\def\Om{\Omega}
\def\E{\mathcal E}

\def\va{\varphi}

\def\O{\Omega}

\def\n{\noindent}

\def\va{\varphi}

\def\pa{\partial}

\def\E{\mathcal E}

\begin{document}
\setlength{\baselineskip}{18truept}
\pagestyle{myheadings}

\markboth{ N. V. Phu and N.Q. Dieu}{Complex $m$-Hessian Type Equations}
\title {Complex $m$-Hessian type equations in weighted energy classes of $m$-subharmonic functions  with given boundary value}
\author{
	Nguyen Van Phu* and Nguyen Quang Dieu**
	\\ *Faculty of Natural Sciences, Electric Power University,\\ Hanoi,Vietnam.\\
	**Department of Mathematics, Hanoi National University of Education,\\ Hanoi, Vietnam;\\
	\\E-mail: phunv@epu.edu.vn\\ and ngquang.dieu@hnue.edu.vn}

\date{}
\maketitle

\renewcommand{\thefootnote}{}

\footnote{2020 \emph{Mathematics Subject Classification}: 32U05, 32W20.}

\footnote{\emph{Key words and phrases}: $m$-subharmonic functions, Complex $m$-Hessian operator, $m$-Hessian type equations, $m$-polar sets, $m$-hyperconvex domain.}

\renewcommand{\thefootnote}{\arabic{footnote}}
\setcounter{footnote}{0}

\begin{abstract}
In this paper, we concern with the existence of solutions of the complex $m-$Hessian type equation $-\chi(u)H_{m}(u)=\mu$ in the class $\mathcal{E}_{m,\chi}(f,\Omega)$  if there exists subsolution in this class, where the given boundary value $f\in\mathcal{N}_m(\Omega)\cap MSH_m(\Om).$
\end{abstract}

\section{Introduction}
 In \cite{Bl1} and \cite{DiKo} the authors introduced $m-$subharmonic functions which are extensions of the plurisubharmonic functions and the complex $m-$Hessian operator $H_m(.) = (dd^c.)^m\wedge \beta^{n-m}$ which is more general than the Monge-Amp\`ere operator $(dd^c.)^n$. In \cite{Ch12}, Chinh introduced the Cegrell classes $\mathcal{F}_m(\Omega)$ and $\mathcal{E}_m(\Omega)$ which are not necessarily locally bounded and  the complex $m-$Hessian operator  is well defined in these classes. In the classes of $m-$subharmonic functions, the complex $m-$Hessian equation $\mu=H_m(u)$ plays important role. Besides solving the $m-$Hessian equation in the case when the measure $\mu$ vanishes on all $m-$polar sets, mathematicians are interested in solving the $m-$Hessian equation when it has subsolution. In \cite{Ch15}, Chinh proved that the complex $m-$Hessian equation has solution in $\mathcal{E}_{m}^{0}(\Omega)$ if it has subsolution in $SH_m(\Omega)\cap L^{\infty}(\Omega).$ Afterthat, in \cite{HP17} the authors proved that the subsolution theorem is true in the class $\mathcal{E}_{m}(\Omega).$ In \cite{Gasmi}, Gasmi extended this result, he  solved complex $m$-Hessian equation in the class $\mathcal{N}_{m}(f)$ if it has subsolution in the class $\mathcal{N}_{m}(\Omega).$ In \cite{AAG20} Amal, Asserda and Gasmi  solved $m-$Hessian type equation $H_{m}(u)=F(u,.)d\mu$ in the class $\mathcal{N}_{m}(f)$ if there exists subsolution in the class $\mathcal{N}_{m}(\Omega)$. Recently, in \cite{PD2023} the authors proved that the complex $m-$Hessian type equations $-\chi(u)H_m(u)=\mu$ has solution in the class $\mathcal{E}_{m,\chi}(\Omega)$ if it has subsolution in this class. Continuing the study in the direction of the  above authors, in this paper, the author will solve complex $m-$Hessian type equation $-\chi(u)H_{m}(u)=\mu$ in the class $\mathcal{E}_{m,\chi}(f,\Omega)$ if it has subsolution in this class where  the given boundary value $f\in\mathcal{N}_m(\Omega)\cap MSH_m(\Om).$ Note that, when $f\equiv 0,$ we get the result in \cite{PD2023}. This results seems to be new even in the plurisubharmonic case. \\
The paper is organized as follows. Besides the introduction,  the paper has other two sections. In Section 2 we recall the definitions and results concerning the $m-$subharmonic functions which were introduced and investigated intensively in recent years by many authors (see \cite{Bl1}, \cite{DiKo}, \cite{SA12}). We also recall the Cegrell classes of $m-$subharmonic functions $\mathcal{F}_m(\Omega)$, $\mathcal{N}_{m}(\Omega)$, $\mathcal{E}_m(\Omega)$ and $\mathcal{E}_{m,\chi}(\Omega)$ which were introduced and studied in \cite{Ch15}, \cite{T19} and \cite{DT23}.   Finally, in Section 3, we solve complex $m-$Hessian type equations $-\chi(u)H_{m}(u)=\mu$ in the class $\mathcal{E}_{m,\chi}(f,\Omega)$ in the case when measure $\mu$ is arbitrary.

\section{Preliminaries}
Throughout this paper, we always denote by $\Om,$ a bounded $m-$hyperconvex domain.
Some elements of the theory of $m$-subharmonic functions and the complex $m$-Hessian operator can be found e.g. in \cite{Bl1}, \cite{DiKo}, \cite{Ch12}, \cite{Ch15}, \cite{SA12} and \cite{T19}. A summary of the properties required for this paper can be found in Preliminaries Section (from subsection 2.1 to subsection 2.6) in \cite{Pjmaa}.\\

\n  We recall some results on weighted $m-$energy classes in \cite{DT23}. Let $\chi: \mathbb{R}^{-}\longrightarrow \mathbb{R}^{-}$ be an increasing function. We put

$$\mathcal{E}_{m,\chi}(\Omega)=\{u\in SH_{m}(\Omega):\exists(u_{j})\in\mathcal{E}_{m}^{0}(\Omega), u_{j}\searrow u, \sup_{j}\int_{\Omega}(-\chi)\circ u_{j}H_{m}(u_{j})<+\infty\}.$$
Note that the weighted $m-$energy classes generalize Cegrell energy  classes $\mathcal{F}_{m,p},\mathcal{F}_{m}.$
\begin{itemize}
	\item When $\chi\equiv -1,$ then $\mathcal{E}_{m,\chi}(\Omega)$ is the class $\mathcal{F}_{m}(\Omega).$\\
	\item When $\chi(t) =-(-t)^{p},$ then $\mathcal{E}_{m,\chi}(\Omega)$ is the class $\mathcal{E}_{m,p}(\Omega).$
\end{itemize}
According Theorem 3.3 in \cite{DT23}, if $\chi\not\equiv 0$ then $\mathcal{E}_{m,\chi}(\Omega)\subset \mathcal{E}_{m}(\Omega)$ which means that the complex $m-$Hessian operator is well - defined on class $\mathcal{E}_{m,\chi}(\Omega)$ and if $\chi(-t)<0$ for all $t>0$ then we have $\mathcal{E}_{m,\chi}(\Omega)\subset \mathcal{N}_{m}(\Omega).$ \\
\n If $\mathcal{K}\in\{\mathcal{E}_m^0(\Omega),\mathcal{E}_m(\Omega),\mathcal{F}_m(\Omega),\mathcal{N}_m(\Omega),\mathcal{F}_m^a(\Omega),\mathcal{N}_m^a(\Omega),\mathcal{E}_{m,\chi}(\Omega \}$ and $f\in\mathcal{E}_m(\Omega).$ We say that a $m-$subharmonic $u$ defined on $\Omega$ belongs to $\mathcal{K}(f)=\mathcal{K}(f,\Omega)$ if there exists a function $\varphi\in\mathcal{K}$ such that $f\geq u\geq \varphi +f.$\\

\n We recall some results that will frequently be used in this paper.

\begin{theorem}[Main Theorem in  \cite{Gasmi}]\label{Gas}
	Let	$\Omega\subset \mathbb{C}^n$ be a bounded $m-$hyperconvex domain and $\mu$ be a positive Borel measure on $\Omega.$ Assume that there exists a function $w\in\mathcal{E}_m(\Omega)$ such that $\mu\leq H_m(w)$ in the sense of currents on $\Omega.$ Then for every function $f\in\mathcal{E}_m(\Omega)\cap MSH_m(\Omega)$, there exists a function $u\in\mathcal{E}_m(\Omega)$ such that $H_m(u)=\mu$ and $f\geq u\geq f+w.$ In particular, if we require $w$ to be in $\mathcal{N}_m(\Omega)$ then $u\in\mathcal{N}_m(f).$
\end{theorem}
\n We recall a version of the comparison principle for a weighted $m-$ Hessian operator (see Theorem 3.8 in \cite{PDjmaa}) .
\begin{theorem}\label{comparison}
	Suppose that the function $t \mapsto \chi (t)$ is nondecreasing in $t.$ 
	Let $u\in\mathcal{N}_{m}(f), v\in\mathcal{E}_{m}(f)$ be such that $-\chi(u)H_{m}(u)\leq -\chi(v)H_{m}(v).$ 
	Assume also that $H_m (u)$ puts no mass on $m$-polar sets.
	Then we have $u\geq v$ on $\Om.$
\end{theorem}
\n We also note the following convergence result for weighted $m-$Hessian operator.
Recall that by Lemma 2.9 in $\cite{T19}$, if a sequence of $m-$subharmonic functions $\{u_j\}$ converges monotonically to a $m-$subharmonic function $u$ then  $u_j\to u$ in $C_m$ as $j\to\infty.$ 
\begin{corollary}[Corollary 3.3 in  \cite{PD2023}]\label{debrecen}
	Let $\chi:\mathbb R^-\longrightarrow\mathbb R^-$ be an increasing continuous function with $\chi(-\infty) >-\infty$. Let $\{u_j, u\}\subset\mathcal E_{m}(\Omega)$, be such that $u_j\geq v$, $\forall j\geq 1$ for some $v\in\mathcal E_{m}(\Omega)$ and that $u_j \to u \in \mathcal E_{m}(\Omega)$ in $C_m$. Then 
	$-\chi (u_j)H_{m}(u_{j})\to -\chi (u)H_{m}(u)$ weakly.
\end{corollary}	
\begin{proposition}[Proposition 2.9 in \cite{Pjmaa}]\label{pro5.2HP17}
	\n  Assume that $u,v,u_k\in\E_m(\Om), k=1,\cdots,m-1$ with $u\geq v $ on $\Om$ and  $T=dd^cu_1\wedge\cdots\wedge dd^cu_{m-1}\wedge\beta^{n-m}.$ Then we have
	$$\ind_{\{u=-\infty\}}dd^cu\wedge T\leq \ind_{\{v=-\infty\}}dd^cv\wedge T.$$
	In particular, if $u,v\in\E_m(\Om)$ are such that $u\geq v$ then for every $m-$polar set $A\subset\Om$ we have
	$$ \int_{A}H_m(u)\leq\int_{A}H_m(v).$$
\end{proposition}
\n We need the following useful approximation result in the class $\mathcal{F}_{m}(f,\Omega).$
\begin{lemma} \label{lm2.6}
	Let  $u \in \mathcal{F}_{m} (f,\Omega)$ with $f\in\mathcal{N}_m(\Omega)\cap MSH_m(\Omega)$ and $\int_{\Omega}H_m(u)<+\infty$. Then there exists a sequence $\{u_j\} \in \mathcal{F}_m (f,\Om)$	such that:	\\
	\n 
	(a) $u_j \downarrow u$ on $\Om;$\\
	\n 
	(b) $H_m (u_j)$ has compact support in $\Om$;\\
	\n 
	(c) $H_m (u_j) \uparrow H_m (u);$\\
	\n (d) $H_m(u_j)$ puts no mass on $m-$polar sets in $\Om.$
\end{lemma}
\begin{proof}
	We fix an element $\rho \in \E^0_m (\Om)\cap C(\Omega)$
	and let $\Omega_j \uparrow \Omega$ be an increasing sequence of relatively compact subsets of $\Omega.$
	For $j \ge 1$ we set
	$\mu_j:= \ind_{\{u>j\rho\} \cap \Om_j} H_m (u).$ Then the measures $\mu_j$ have the following properties:
	
	\n 
	(i) $\mu_j$ has compact support in $\Om;$
	
	\n 
	(ii) $\mu_j \le \mu_{j+1} \le H_m (u);$ 
	
	\n 
	(iii) $\mu_j$ puts no mass on $m-$polar sets in $\Om$ (by Lemma 2.16 in \cite{Gasmi}); 
	
	\n 
	(iv) $\int\limits_\Om d\mu_j \le \int\limits_{\{u>j\rho\}} H_m(u) \le j^m \int\limits_{\Om} H_m (\rho) <\infty$ 
	(by Lemma 5.5 in $\cite{T19}$).
	
	\n 
	It follows from the hypothesis $u \in \mathcal{F}_{m} (f,\Omega)$ and $f\in\mathcal{N}_m(\Omega)$ that $u\in\mathcal{N}_m(\Omega).$ Thus using (ii), (iv) and the main Theorem in $\cite{Gasmi}$ we can find $u_j \in \mathcal{N}_m (f,\Om)$ such that $H_m (u_j)=\mu_j.$ We have $\int_{\Omega}H_m(u_j)=\int_{\Omega}d\mu_j<\infty.$ Thus, by Theorem 3.1 in \cite{Pjmaa}  $u_j \in \mathcal{F}_m (f,\Om).$  Obviously, we have  $u_j$ satifies properties (b), (c) and (d).\\ So it remains to prove that $u_j \downarrow u$ on $\Om.$
	Indeed, by Theorem 3.8 in $\cite{PDjmaa}$, (ii) and (iii) we get $u_j \geq u_{j+1}\geq u.$ 
	Set $v:= \lim\limits_{j \to \infty} u_j$ then we have $v \ge u.$ Moreover, by Corollary \ref{debrecen} we deduce that $H_m(u_j)\to H_m(v)$ weakly as $j\to\infty.$ Coupling with 
	the construction of $u_j$ we have $H_m (v)=H_m (u).$ By the Theorem 2.10 in \cite{Gasmi} (see Theorem 3.6 in \cite{ACCH} for the case of plurisubharmonic functions) we obtain $u=v$, and so we have $u_j \downarrow u$ on $\Om.$ 
	The proof is completed.	
\end{proof}
\n In connection to Lemma \ref{lm2.6} we present the following result which might be of independent interest.
This result was also used implicitly in the proof of \cite{PD2023}.
\begin{lemma} \label{lm2.5}
	Let  $u \in \mathcal F_{m} (f,\Omega), f\in \mathcal{E}_m(\Omega)\cap MSH_m(\Omega).$  Assume that the support of $H_m (u)$ is a compact subset of $\Om,$	then there exist an open subset $\Om' \Subset \Om$
	a sequence $\{u_j\} \in \mathcal{E}_m^0 (f,\Om)$ having the following properties:	\\
	\n 
	(a) $u_j \downarrow u$ on $\Om;$\\
	\n 
	(b) $H_m (u_j)$ has compact support in $\overline{\Om'}$;\\
	\n 
	(c) $H_m (u_j)$ puts no mass on $m-$polar subsets of $\Om.$ 
\end{lemma}
\begin{proof}
	Choose a domain $\Om' \Subset \Om$ such that $\Om'$ contains the support of $H_m (u), \pa \Om'$ is $\mathcal C^1-$smooth
	and $H_m (u)$ puts no mass on $\pa \Om'.$  According to Proposition 2.12 in \cite{Gasmi}, we may find a sequence $v_j \in \E^0_m (\Om, f)$ such that $v_j \downarrow u$ on $\Om$.
	Set
	$$u_j:=  \sup \{\va: \va \in SH_m^{-} (\Om), \va|_{\Om'} \le v_j, \va \le f\}.$$
	By maximality of $f$ we see that $H_m (u_j)=0$ on $\Om \setminus \overline{\Om'}.$
	Since $u_j \ge v_j$ on $\Om$ there exists a function $\xi_j \in\mathcal{E}_m^0(\Om)$ such that $u_j\geq f+\xi_j.$ According to Proposition $\ref{pro5.2HP17}$ with note that $f\in MSH_m(\Om),$ for all $m$-polar set $A\subset\Om$ we also have
	$$\int_{A}H_m(u_j)\leq\int_{A}H_m(f+\xi_j)= 0,$$
	where the last inequality is due to Lemma 5.6 in \cite{HP17}.
	That means $H_m(u_j)$ vanishes on all $m$-polar sets. So we are done.
	
	For $(a)$, we first observe that
	$u_j \downarrow:=v \ge u$ on $\Om.$ Moreover, since $u_j \ge v_j$ on $\Om$ we infer that $u_j=v_j$ on $\Om'.$
	Thus $v= u$ on $\Om'$, and so $v=u$ on $\partial \Om'$.
	Now we define
	$$\tilde u:= (\sup \{\va \in SH^{-}_m (\Om'): \va^*|_{\pa \Om'} \le u\})^*.$$
	Since $\pa \Om'$ is $\mathcal C^1-$smooth we have $\tilde u \in SH_m (\Om')$ and $\tilde u \ge u$ on $\Om'.$
	Hence the function 
	$$\hat u:= \begin{cases} u & \text{on}\  \Om \setminus \Om'\\
	\tilde u & \text{on}\ \Om'
	\end{cases}$$
	belongs to $SH_m (\Om),$ and since $\hat u \ge u$ on $\Om,$ we infer that $\hat u \in \E_m (\Om).$ Observe also that
	$H_m (\hat u)$ is supported on $\pa \Om'.$
	Observe that by Proposition $\ref{pro5.2HP17}$ and the choice of $\Om'$ we have 
	\begin{equation} \label{eq000}
	\int\limits_{\{\hat u=-\infty\} \cap \partial \Om'} H_m (\hat u) \le \int\limits_{\{\hat u=-\infty\} \cap \partial \Om'} H_m (u)=0.
	\end{equation}
	Since  $\{v >\hat u >-\infty\} \cap \pa \Om'$ is empty, in view of (\ref{eq000}), we may apply Lemma 3.1 in \cite{Gasmi}
	to conclude that $\hat u \ge v$ on $\Om.$ So in particular $u \ge v$ on $\Om \setminus \Om'.$
	Therefore $u=v$ on $\Om.$
	Thus we obtain $u_j \downarrow v=u$ on $\Om.$ That completes the proof of our lemma.
\end{proof}
\section{ Complex $m-$Hessian equations in the class $\mathcal{E}_{m,\chi}(f,\Omega)$}
In this section, we assume that  $\chi:  \mathbb{R}^{-}\to \mathbb{R}^{-} $ is a  nondecreasing continuous function  such that  $\chi(t)<0$ for all  $t<0$. 
We first concern with the complex $m-$Hessian equations $-\chi(u)H_{m}(u)=\mu$ in the class $\mathcal{E}_{m,\chi}(f,\Omega)$ when $\mu$ puts no mass on $m-$polar sets. 
\begin{theorem}\label{mpolar}
	Let $\mu$ be a nonnegative, finite measure which puts no mass on $m-$polar sets. Then the complex $m-$Hessian type equation $-\chi(u)H_{m}(u)=\mu$ has solution in the class $\mathcal{E}_{m,\chi}(f,\Omega),$ where $f\in \mathcal{E}_m(\Omega)\cap MSH_m(\Omega).$
\end{theorem}
\begin{remark} According to Theorem 3.1 in \cite{AAG20}, if $\mu$ is a nonnegative measure which puts no mass on $m-$polar sets and the complex equation $-\chi(u)H_{m}(u)=\mu$ has a subsolution in $\mathcal{N}_m^a(\Omega)$ then it has a solution in $\mathcal{N}_{m}(f,\Omega).$  Main Theorem in \cite{HQ21} also proved that if $\mu$ be a nonnegative measure which puts no mass on $m-$polar sets and if the equation $-\chi(u)H_{m}(u)=\mu$ has a subsolution in $\mathcal{E}_m(\Omega)$ then it has solution which belongs to $\mathcal{E}_{m}(\Omega).$  Our theorem \ref{mpolar}  does not require the existence of a subsolution, but instead
	we need finiteness of the measure $\mu.$  On the other hand,  the solution we found is somewhat more precise since it is contained in $\mathcal{E}_{m,\chi}(f,\Omega)\subset \mathcal{N}_m(f,\Omega)\subset\mathcal{E}_m(\Omega).$ 
\end{remark}
\begin{proof}
	By Theorem 5.3 in \cite{Ch15} we can find $\varphi\in\mathcal{E}_{m}^{0}(\Omega)$ and $0\leq h\in L^{1}_{loc}(H_{m}(\varphi))$ such that $\mu=hH_{m}(\varphi).$
Set $\mu_{j}=1_{\Omega_{j}}\min(h,j)H_{m}(\varphi)$ where $\{\Omega_{j}\}$ be a fundamental sequence of $\Omega.$ \\
Choose nondecreasing functions
$\chi_{j}\in C^{\infty}(\mathbb{R}^{-})$ such that $-\chi_{j}\searrow -\chi.$ Put $\gamma(t) = \frac{1}{\chi(t)}$ and $\gamma_{j}(t)=\dfrac{1}{\chi_{j}(t)}$. We have that $\gamma(t)$ is a nonincreasing function and nonincreasing functions $\gamma_{j}\in C^{\infty}(\mathbb{R}^{-})$ satisfying -$\gamma_{j}(t)\nearrow -\gamma(t) .$    \\ 
Note that $-\gamma_j$ is above bounded on $\Omega_j$ so using 
Proposition 3.4 in \cite{AAG20}, we can find $u_{j}\in\mathcal{N}_{m}(f)$ such that 
$$H_{m}(u_{j})=-\gamma_j(u_j)d\mu_{j}=\dfrac{d\mu_{j}}{-\chi_j(u_j)}.$$
It follows that $$-\chi_j(u_{j})H_{m}(u_{j})= \mu_{j}.$$
Therefore, we have
$$-\chi_j(u_{j})H_{m}(u_{j})= \mu_{j}\leq \mu_{j+1}=-\chi_{j+1}(u_{j+1})H_{m}(u_{j+1})\leq -\chi_{j}(u_{j+1})H_{m}(u_{j+1}).$$
By Theorem \ref{comparison} we have $u_{j}\searrow u.$ We will prove that $u\in\mathcal{E}_{m,\chi}(f,\Omega)$ which satisfies $$-\chi(u)H_{m}(u)=\mu.$$
Firstly, we  prove that $u_{j}\in\mathcal{E}_{m}^{0}(f).$ 
Applying Proposition 3.4 in \cite{AAG20} one again (in the case $f\equiv 0$), we can find $\varphi_{j}\in\mathcal{F}_{m}^{a}(\Omega)$ such that 
$$H_{m}(\varphi_{j})=-\gamma_j(\varphi_j) d\mu_{j}=-\dfrac{d\mu_{j}}{\chi_j(\varphi_j)}.$$
This implies that $$-\chi_j(\varphi_{j})H_{m}(\varphi_{j})= \mu_{j}.$$
On the other hand, since $\Omega_j\Subset\Omega$, we obtain $\varphi_j\leq A(j)<0$ on $\Omega_j.$ Note that $-\gamma_j$ is a nondecreasing continuous function, we deduce that $-\gamma_j(\varphi_j)\leq-\gamma_j(A(j))\leq -\gamma(A(j))=B(j)$ on $\Omega_j.$ Hence
$$H_{m}(\varphi_{j})=-\gamma_j(\varphi_{j}) \mu_{j}\leq jB(j)H_{m}(\varphi)
=H_m(\sqrt[m]{jB(j)}\varphi).$$
It follows  Theorem \ref{comparison}  that $\varphi_{j}\geq \sqrt[m]{jB(j)}\varphi.$ Since $\varphi\in\mathcal{E}_m^0(\Omega)$ we obtain $\varphi_j\in\mathcal{E}_{m}^{0}(\Omega).$ 
Moreover, we have
$$-\chi_j(u_{j})H_{m}(u_{j})=\mu_{j}=-\chi_j(\varphi_{j})H_{m}(\varphi_{j})\leq -\chi_j(f+\varphi_{j})H_{m}(f+\varphi_{j}) $$
and $u_{j}, f+\varphi_{j}\in\mathcal{N}_{m}^{a}(f)$ then Theorem $\ref{comparison}$ implies that $u_{j}\geq f+\varphi_{j}.$ So we have $u_{j}\in\mathcal{E}_{m}^{0}(f)$ as the desired.\\
Secondly, we prove that $u\in\mathcal{E}_{m,\chi}(f).$ Indeed, we have 
$$-\chi_j(\varphi_{j})H_{m}(\varphi_{j})=d\mu_{j}\leq d\mu_{j+1} =-\chi_{j+1}(\varphi_{j+1})H_{m}(\varphi_{j+1})\leq -\chi_{j}(\varphi_{j+1})H_{m}(\varphi_{j+1}).$$ 
According to Theorem \ref{comparison} we see that $\{\varphi_{j}\}$ is decreasing  and we assume that $\psi=\lim\limits_{j\to\infty}\varphi_{j}.$ Note that $\varphi_j\in\mathcal{E}_m^0(\Omega).$ Moreover, we have 
$$\sup_{j\geq 1}\int\limits_{\Omega}-\chi(\varphi_j)H_m(\varphi_{j})\leq \sup_{j\geq 1}\int\limits_{\Omega}-\chi_j(\varphi_j)H_m(\varphi_{j}) =\sup_{j\geq 1}\int\limits_{\Omega}d\mu_{j}\leq\mu(\Omega)<\infty.$$ Therefore, we obtain $\psi\in\mathcal{E}_{m,\chi}(\Omega).$ 
It follows from $f\geq u_{j}\geq f+\varphi_{j}$ that $f\geq u\geq f+\psi$ and we get $u\in\mathcal{E}_{m,\chi}(f)$ as desired. \\
Thirdly, we prove that $-\chi(u)H_{m}(u)=\mu.$
Indeed, we have $-\chi(u_j)H_m(u_j)=\mu_j.$ Repeating the argument as in the last part in the proof of Theorem 4.1 in \cite{PD2023} we have $\lim\limits_{j\to\infty}H_m(u_j)=-\gamma(u) \mu.$
On the other hand, since $u_j\searrow u\in \mathcal{E}_{m,\chi}(f)\subset\mathcal{E}_m(\Omega),$ according to Theorem 3.8 in \cite{HP17} we obtain $H_m(u_j)$ converges weakly to $H_m(u)$ as $j\to\infty.$ So we have
$$H_m(u)=-\eta(u)\mu \Rightarrow -\chi(u)H_m(u)=\mu.$$
The proof is complete. 
\end{proof}
\n The next result deals with the case $\mu$ is a arbitrary measure with finite total mass. 
\begin{theorem}
\label{theorem4.2}
Let  $\mu$ be a non-negative finite measure on $\Omega.$ Assume that there exists a function $w\in\mathcal E_{m,\chi}(f,\Omega)$ with $\mu\leq -\chi(w)H_{m}(w)$, where the given boundary $f\in\mathcal{N}_m(\Omega)\cap MSH_m(\Omega)$. Then there exists a function $u\in\mathcal{E}_{m,\chi}(f,\Omega)$ such that $u\geq w$ and $-\chi(u)H_{m}(u)=\mu$.
\end{theorem}
\begin{remark}
 According to Main Theorem in \cite{AAG20}, we only achieve a solution $u\in\mathcal{N}_m(f)$ and $u\geq f+w.$ In Theorem
 \ref{theorem4.2} we have finer information about this solution $u.$		
\end{remark}
\begin{proof}
	 We consider two cases.
	
	\n {\em Case 1.} Assume that $\chi(-\infty)>-\infty$. 
Using Theorem 2.15 in \cite{Gasmi} we may decompose   $\mu = \alpha +\nu$, where $\alpha$ and $\nu$ are Radon measures defined on $\Omega$  such that $\alpha$ vanishes on all $m-$polar sets and $\nu$ is carried by an $m-$polar set. It follows from $w\in\mathcal{E}_{m,\chi}(f,\Omega)$ and $f\in\mathcal{N}_m(\Omega)$ that  $w\in\mathcal{N}_m(\Omega)$. The  Theorem  \ref{Gas} implies that there exists $v\in\mathcal N_m(f,\O)$ such that $v\geq f+ w$ and $\nu=H_m(v).$ 
Note that, $$\int_{\Omega}H_m(v)=\int_{\Omega}d\nu\leq \int_{\Omega}d\mu<+\infty.$$ Thus, by Theorem 3.1 in \cite{Pjmaa} we infer that $v\in\mathcal F_m(f,\O).$ According to Lemma \ref{lm2.6}, there exist $v_j\in\mathcal{F}_m(f,\Omega)$ such that $v_j\searrow v,$ supp$H_m(v_j)\Subset\Omega, H_m(v_j)$ puts no mass on  $m-$polar sets and $\sup\limits_{j \ge 1} \int_{\Om} H_m (v_j)<\infty.$

Using Theorem \ref{mpolar} we can find  $u_j\in\mathcal{E}_ {m, \chi}(f,\Omega)$ that satisfies
	\begin{equation}\label{eqp}
	-\chi (u_j) H_m(u_j) = \alpha + H_m(v_j).
	\end{equation}
Observe that for $j \ge 1$ we have	
$$-\chi (u_j) H_m(u_j) \le -\chi (u_{j+1}) H_m(u_{+1}) \le -\chi (w) H_m(w).$$
It then follows from Theorem \ref{comparison} that $u_j\searrow u\geq w\in\mathcal{E}_{m,\chi}(f,\Omega).$ This implies that $u\in\mathcal{E}_{m,\chi}(f,\Omega)\subset\mathcal{E}_m(\Omega).$ Note that we have $ \chi(-\infty)>-\infty$ , by letting $j \to \infty$ in (\ref{eqp}) and using Corollary \ref{debrecen}, we obtain
$$-\chi (u) H_m(u) = \alpha + H_m(v)=\mu.$$
\n {\em Case 2.} Assume that $\chi(-\infty)=-\infty$. It follows from the hypothesis $w\in\mathcal E_{m,\chi}(f,\Omega)$ that there exists a function $\psi\in\mathcal{E}_{m,\chi}(\Omega)$ such that $f\geq w\geq f+\psi.$ By Theorem 3.7  in \cite{DT23} we have  $\psi\in\mathcal E^a_{m}(\Omega).$ Note that we have $f\in MSH_m(\Omega)$ so by Theorem 1.2 in \cite{Bl1} we obtain $H_m(f)=0.$ Therefore, for every $m-$polar set $A\subset\Omega,$ by Proposition \ref{pro5.2HP17} and Lemma 5.6 in \cite{HP17} we infer that
$$\int_{A}H_m(w)\leq\int_{A}H_m(f+\psi)=0.$$
This means that $H_m(w)$ vanishes on pluripolar sets and so is $\mu$. Thus, by Theorem \ref{mpolar} there  exists a function $u\in\mathcal E_{m,\chi}(f,\Omega)$ such that $ -\chi(u)H_{m}(u)=\mu$. The proof is complete.
\end{proof}

\section*{Declarations}
\subsection*{Ethical Approval}
This declaration is not applicable.
\subsection*{Competing interests}
The authors have no conflicts of interest to declare that are relevant to the content of this article.
\subsection*{Authors' contributions }
Nguyen Van Phu and Nguyen Quang Dieu together studied  the manuscript.
\subsection*{Availability of data and materials}
This declaration is not applicable.

\end{document}